\documentclass{article}
\usepackage{latexsym,amsthm,amssymb,amscd,amsmath}
\newtheorem{thm}{Theorem}[section]
\newtheorem{cor}[thm]{Corollary}

\newtheorem{lem}[thm]{Lemma}

\newcommand{\be}{\begin{equation}}
\newcommand{\ee}{\end{equation}}
\newcommand{\ben}{\begin{enumerate}}
\newcommand{\een}{\end{enumerate}}
\newcommand{\beq}{\begin{eqnarray}}
\newcommand{\eeq}{\end{eqnarray}}
\newcommand{\beqn}{\begin{eqnarray*}}
\newcommand{\eeqn}{\end{eqnarray*}}

\newcommand{\pa}{\partial}

\newcommand{\pxi}{ {\pa \over \pa x^i}}

\title{On $m$-th Root Finsler Metrics}
\author{A. Tayebi and B. Najafi}
\begin{document}

\maketitle
\begin{abstract}
In this paper, we characterize locally dually flat and Antonelli $m$-th root Finsler metrics. Then, we show
that every  $m$-th root Finsler metric of isotropic mean Berwald curvature  reduces to a weakly Berwald metric.\\\\
{\bf {Keywords}}: Antonelli metric, locally dually flat metric, isotropic mean Berwald metric.\footnote{ 2000 Mathematics subject Classification: 53C60, 53C25.}
\end{abstract}

\section{Introduction}
The theory of $m$-th root metrics has been developed by H. Shimada \cite{Shim}, and applied to Biology as an ecological metric \cite{AIM}. It is regarded as a direct generalization of Riemannian metric in the sense that the second root metric is a Riemannian metric.  The third and fourth root metrics are called the cubic metric and quartic metric, respectively.

Recently studies  show that the theory of  $m$-th root Finsler metrics plays a very important role in physics, theory of space-time structure, gravitation, general relativity and seismic ray theory \cite{Asan}\cite{HD}\cite{MatShim}\cite{Pav1}\cite{Pav2}. For quartic metrics, a study of the geodesics and of the related geometrical objects is made by S. Lebedev \cite{Leb}. Also, Einstein equations for some relativistic models relying on such metrics are studied by V. Balan and N. Brinzei in two papers \cite{Mangalia}, \cite{Cairo1}. Tensorial
connections for such spaces have been recently studied by L. Tamassy \cite{Tamassy}. B. Li and Z. Shen study locally projectively flat fourth root metrics under some irreducibility condition \cite{shLi1}. Y. Yu and Y. You show that an $m$-th root Einstein Finsler metric is Ricci-flat \cite{YY}.

Let $(M, F)$ be a Finsler manifold of dimension $n$, $TM$ its tangent bundle and $(x^{i},y^{i})$ the
coordinates in a local chart on $TM$. Let $F$ be a scalar function on $%
TM$ defined by $F=\sqrt[m]{A}$, where $A$ is given by
\begin{equation}
A:=a_{i_{1}\dots i_{m}}(x)y^{i_{1}}y^{i_{2}}\dots y^{i_{m}},
\label{1.1}
\end{equation}
with $a_{i_{1}\dots i_{m}}$ symmetric in all its indices \cite{Shim}. Then $F$ is called an $m$-th root Finsler metric.

Let  $F$ be an $m$-th root Finsler metric on an open subset $U\subset \mathbb{R}^n$. Put
\[
A_i=\frac{\pa A}{\pa y^i}, \ \ \textrm{ and} \ \  A_{ij}=\frac{\pa^2 A}{\pa y^i\pa y^j}.
\]
Suppose that the matrix $(A_{ij})$ defines a positive definite tensor and $(A^{ij})$ denotes its inverse.   Then the following  hold
\begin{eqnarray}
&&g_{ij}= \frac{A^{\frac{2}{m}-2}}{m^2}[mAA_{ij}+(2-m)A_iA_j],\label{g}\\
&&g^{ij}= A^{-\frac{2}{m}}[mAA^{ij}+\frac{m-2}{m-1}y^iy^j],\\
&&y^iA_i=mA, \ \ y^iA_{ij}=(m-1)A_j,\ \ y_i=\frac{1}{m}A^{\frac{2}{m}-1}A_i,\\
&&A^{ij}A_{jk}=\delta^i_k,\ \ A^{ij}A_i=\frac{1}{m-1}y^j, \ \ A_iA_jA^{ij}=\frac{m}{m-1}A,\\
&&A_0:=A_{x^m}y^m, \ \ A_{0l}:=A_{x^my^l}y^m.
\end{eqnarray}
A Finsler metric $F=F(x,y)$ on a manifold $M$ is said to be locally dually flat, if at any point there is a standard coordinate system
 $(x^i,y^i)$ in $TM$ such that $[F^2]_{x^ky^l}y^k=2[F^2]_{x^l}$. In this case, the coordinate $(x^i)$ is called an adapted local coordinate system \cite{shen1}. Here, we characterize locally dually flat $m$-th root Finsler metrics.
\begin{thm} \label{main thm 1}
Let  $F$ be an $m$-th root Finsler metric on an open subset $U\subset \mathbb{R}^n$. Then,  $F$ is a locally dually flat metric if and only if the following holds
\be
A_{x^l}=\frac{1}{2A}\Big\{(\frac{2}{m}-1)A_lA_0+AA_{0l}\Big\}.\label{d0}
\ee
Moreover, suppose that $A$ is irreducible. Then $F$ is locally dually flat if and only if there exists a 1-form $\theta=\theta_l(x)y^l$
on $U$ such that the following holds
\be
A_{x^l}=\frac{1}{3m}\Big\{2\theta A_l+mA\theta_l\Big\}.\label{d10}
\ee
\end{thm}
A Finsler metric $F$ is called an Antonelli metric, if  there is a local coordinate system in $TM$, in which the spray coefficients are dependent on direction alone. In this case, the spray coefficients of $F$ are given by $G^i=\frac{1}{2}\Gamma^i_{jk}(y) y^jy^k$.

Antonelli metrics  were first introduced by P. L. Antonelli for some studies in Biology and Ecology \cite{AIM}. Antonelli calls them $y$-Berwald metrics. This class of  metrics arises in time sequencing change models in the evolution of colonial systems. Here, we characterize  $m$-th root Antonelli metrics. More precisely, we prove the following.
\begin{thm} \label{main thm 2}
Let  $F$ be an $m$-th root Finsler metric on an open subset $U\subset \mathbb{R}^n$. Then, $F$ is an Antonelli metric if and only if there exist
functions $\Gamma^i_{lk}$ dependent only on direction such that the following holds
\be
A_{x^l}=[\Gamma^i_{lk}y^k+\frac{1}{2}\Gamma^i_{jk,l}y^jy^k] A_{i},
\ee
where $\Gamma^i_{jk,l}=\frac{\pa \Gamma^i_{jk} }{\pa y^l}$.
\end{thm}
A Finsler metric $F$ is  called a Berwald metric if  $G^i$  are quadratic in $y\in T_xM$  for any $x\in M$ or  equivalently  Berwald curvature  vanishes. In this case, we have $G^i = {1\over 2} \Gamma^i_{jk}(x)y^jy^k$. Hence, Berwald metrics can be considered as $x$-Berwald metrics. The $E$-curvature is defined by the trace  of the Berwald curvature. A Finsler metric $F$ is called of isotropic mean Berwald curvature if ${\bf E}= \frac{n+1}{2} c F^{-1}{\bf h}$, where $c=c(x)$ is a scalar function on $M$ and ${\bf h}$ is the angular metric. If $c=0$, then $F$ is called weakly Berwald metric.
\begin{thm} \label{main thm 3}
Let  $F$ be an $m$-th root Finsler metric on an open subset $U\subset \mathbb{R}^n$ with $n\geq 2$. Suppose that $F$ is of isotropic mean Berwald curvature. Then $F$ is a weakly Berwald metric.
\end{thm}

\section{Preliminaries}
Let $M$ be an $n$-dimensional $ C^\infty$ manifold. Denote by $T_x M $
the tangent space at $x \in M$, by $TM=\cup _{x \in M} T_x M $
the tangent bundle of $M$ and  by  $TM_{0} = TM \setminus \{ 0 \}$ the slit tangent bundle of $M$.

A Finsler metric on a manifold $M$ is a function $ F:TM\rightarrow [0,\infty)$ with the following properties:\\
(i) $F$ is $C^\infty$ on $TM_{0}$;\\
(ii) $F(x,\lambda y)=\lambda F(x,y) \,\,\, \forall \lambda>0$,\ \ $y\in TM$;\\
(iii) for each $y\in T_xM$, the following quadratic form $g_y$ on $T_xM$  is positive definite,
\[
g_{y}(u,v):={1 \over 2} \left[  F^2 (y+su+tv)\right]|_{s,t=0}, \ \
u,v\in T_xM.
\]

\bigskip

In  local coordinates $(x^i,y^i)$, the vector filed ${\bf G}=y^i\frac{\pa}{\pa x^i}-2G^i\frac{\pa}{\pa y^i}$ is a global vector field on $TM_0$, where $G^i=G^i(x,y)$ are local functions on $TM_0$  given by
\be \label{spray}
G^i:=\frac{1}{4}g^{il}\Big\{ \frac{\partial^2F^2}{\partial x^k \partial y^l}y^k-\frac{\partial F^2}{\partial x^l}\Big\},\ \ y\in T_xM.
\ee
The vector field ${\bf G}$ is called the associated spray to $(M,F)$ \cite{TP}.

\bigskip

For $y \in T_xM_0$, define ${\bf B}_y:T_xM\otimes T_xM \otimes T_xM\rightarrow T_xM$ and ${\bf E}_y:T_xM \otimes T_xM\rightarrow \mathbb{R}$ by ${\bf B}_y(u, v, w):=B^i_{\ jkl}(y)u^jv^kw^l{{\partial } \over {\partial x^i}}|_x$ and ${\bf E}_y(u,v):=E_{jk}(y)u^jv^k$
where
\[
B^i_{\ jkl}:={{\partial^3 G^i} \over {\partial y^j \partial y^k \partial y^l}},\ \ \ E_{jk}:={{1}\over{2}}B^m_{\ jkm}.
\]
$\bf B$ and $\bf E$ are called  Berwald curvature and mean Berwald curvature, respectively.  Then $F$ is called a Berwald metric and weakly Berwald metric if ${\bf{B}}=0$ and ${\bf{E}}=0$, respectively \cite{TP}.

 A Finsler metric $F$ on an $n$-dimensional manifold $M$ is said to be isotropic mean Berwald metric or of isotropic $E$-curvature if
\[
E_{ij} = \frac{n+1}{2} cF^{-1} h_{ij},
\]
where $h_{ij}= g_{ij} - g_{ip}y^p g_{jq}y^q$ is the angular metric.
\section{Proof of Theorem \ref{main thm 1}}
Information geometry has emerged from investigating the geometrical structure of a family of probability distributions and has been applied successfully to various areas including statistical inference, control system theory and multi-terminal information theory \cite{am}\cite{amna}. Dually flat Finsler metrics form a special and valuable class of Finsler metrics in Finsler information geometry, which plays a very   important role in studying flat Finsler information structure \cite{shen1}.

A Finsler metric $F=F(x,y)$ on a manifold $M$ is said to be locally dually flat, if at any point there is a standard coordinate system
 $(x^i,y^i)$ in $TM$ such that $L=F^2$ satisfies
\[
L_{x^ky^l}y^k=2L_{x^l}.
\]
In this case, the coordinate $(x^i)$ is called an adapted local coordinate system. Every locally Minkowskian metric $F$ satisfies trivially the above equation, hence $F$ is locally dually flat.

In order to find explicit examples of locally dually flat metrics, we consider $m$-th root Finsler metrics.

\bigskip

\noindent{\bf{Proof of Theorem \ref{main thm 1}}}: Let $F$ be a locally dually flat metric
\be
(A^{\frac{2}{m}})_{x^ky^l}y^k=2(A^{\frac{2}{m}})_{x^l}.\label{d1}
\ee
We have
\begin{eqnarray}
(A^{\frac{2}{m}})_{x^l}\!\!\!\!&=&\!\!\!\!\!\! \frac{2}{m}A^{\frac{2-m}{m}}A_{x^l},\label{d2}\\
(A^{\frac{2}{m}})_{x^ky^l}\!\!\!\!&=&\!\!\!\!\!\! \frac{2}{m}A^{\frac{2-2m}{m}}\Big[(\frac{2}{m}-1)A_lA_{x^k}+AA_{x^ky^l}\Big].\label{d3}
\end{eqnarray}
By (\ref{d1}), (\ref{d2}) and (\ref{d3}), we have (\ref{d0}). The converse is trivial.

Now, suppose that $A$ is irreducible. One can rewrite (\ref{d0}) as follows
\be
A(2A_{x^l}-A_{0l})=(\frac{2}{m}-1)A_lA_0.\label{d11}
\ee
Irreducibility of $A$ and $deg(A_l)=m-1$ imply that there exists a 1-form $\theta=\theta_l y^l$ on $U$ such that
\be
A_0=\theta A.\label{d12}
\ee
Plugging (\ref{d12}) into (\ref{d0}), we get
\be
A_{0l}=A\theta_l+\theta A_l-A_{x^l}.\label{d13}
\ee
Substituting (\ref{d12}) and (\ref{d13}) into (\ref{d0}) yields (\ref{d10}).  The converse is a direct computation. This
completes the proof.
\qed

\bigskip

\begin{cor}
Let  $F=\sqrt{A}$ be a Riemannian metric on open subset $U\subset \mathbb{R}^n$. Then,  $F$ is a locally dually flat metric if and only if  there exists a 1-form $\theta$
on $U$ such that the following holds
\be
A_{x^l}=\frac{1}{3}(\theta A)_{y^l}.\label{d14}
\ee
In this case, the spray coefficients of $F$ are given by
\be \label{spary2}
G^i=\frac{1}{12}\theta^iF^2+\frac{1}{6}\theta y^i,
\ee
where $\theta^i=2A^{ik}\theta_k$.
\end{cor}
\begin{proof} Putting $m=2$ in (\ref{d10}), we get (\ref{d14}). Plugging $A=a_{ij}y^iy^j$ into (\ref{d14}), one can obtain
\be
3\frac{\pa a_{ij}}{\pa x^l}=\theta_la_{ij}+\theta_ia_{lj}+\theta_ja_{il}. \label{d15}
\ee
Substituting (\ref{d15}) into (\ref{spray}), we get (\ref{spary2}).
\end{proof}

\section{Proof of Theorem \ref{main thm 2}}
In this section, we deal with Antonelli metrics. First, we remark the following.
\begin{lem}\label{YY}{\rm (\cite{YY})}
\emph{Let  $F$ be an $m$-th root Finsler metric on an open subset $U\subset \mathbb{R}^n$. Then the spray coefficients of $F$ are given by}
\be
G^i=\frac{1}{2}(A_{0j}-A_{x^j})A^{ij}.\label{YY}
\ee
\end{lem}

Now, we are going to prove  Theorem \ref{main thm 2}. First, we prove the following.

\begin{lem}\label{Lem1}
Let  $F$ be an $m$-th root Finsler metric on an open subset $U\subset \mathbb{R}^n$. If $F$ is an Antonelli metric, then the following holds
\be
A_{x^l}=[\Gamma^i_{lk}y^k+\frac{1}{2}\Gamma^i_{jk,l}y^jy^k] A_{i},\label{a0}
\ee
where $\Gamma^i_{jk,l}=\frac{\pa \Gamma^i_{jk} }{\pa y^l}$.
\end{lem}
\begin{proof}
Let $F$ be an Antonelli metric metric. Then
\be
G^i=\frac{1}{2}\Gamma^i_{jk}(y) y^jy^k.\label{a1}
\ee
Plugging (\ref{a1}) into (\ref{YY}), we get
\be
\Gamma^i_{jk}y^jy^k=(A_{0j}-A_{x^j})A^{ij}.\label{a2}
\ee
Multiplying (\ref{a1})  with $A_{il}$ and $A_{i}$ implies that
\begin{eqnarray}
\Gamma^i_{jk}y^jy^k A_{il}\!\!\!\!&=&\!\!\!\!\ A_{0l}-A_{x^l},\label{a3}\\
\Gamma^i_{jk}y^jy^k A_{i}\!\!\!\!&=&\!\!\!\!\ A_0.\label{a4}
\end{eqnarray}
Differentiating of (\ref{a4}) with respect to $y^l$, we have
\be
[\Gamma^i_{jk,l}y^jy^k+2\Gamma^i_{lk}y^k] A_{i}+\Gamma^i_{jk}y^jy^k A_{il}=A_{0l}+A_{x^l}.\label{a5}
\ee
Subtracting (\ref{a5}) from (\ref{a3}) yields
\be
A_{x^l}=[\Gamma^i_{lk}y^k+\frac{1}{2}\Gamma^i_{jk,l}y^jy^k] A_{i}.\label{a6}
\ee
Then we get the proof.
\end{proof}
\begin{lem}\label{Lem2}
Let  $F$ be an $m$-th root Finsler metric on an open subset $U\subset \mathbb{R}^n$. Suppose that $F$  satisfies the following equation
\be
A_{x^l}=[\Gamma^i_{lk}y^k+\frac{1}{2}\Gamma^i_{jk,l}y^jy^k] A_{i},\label{a00}
\ee
where $\Gamma^i_{lk}$ are functions only of direction. Then, $F$ is an Antonelli metric.
\end{lem}
\begin{proof}
Now suppose that (\ref{a00}) holds. By differentiating of (\ref{a00}) with respect to $y^h$, we have
\begin{eqnarray}
\nonumber A_{x^ly^h}=[\Gamma^i_{lk,h}y^k+\Gamma^i_{lh}\!\!\!\!&+&\!\!\!\! \frac{1}{2}\Gamma^i_{jk,l,h}y^jy^k+\Gamma^i_{hk,l}y^k] A_{i}\\
\!\!\!\!&+&\!\!\!\! [\frac{1}{2}\Gamma^i_{jk,l}y^jy^k+\Gamma^i_{lk}y^k]A_{ih}.\label{a7}
\end{eqnarray}
Contracting (\ref{a7}) with $y^l$ and 0-homogeneity of functions $\Gamma^i_{lk}$ yield
\be
A_{0h}=[\frac{1}{2}\Gamma^i_{lk,h}y^ly^k+\Gamma^i_{lh}y^l]A_{i}+\Gamma^i_{lk}y^ky^lA_{ih}.\label{a8}
\ee
Substituting (\ref{a00}) into (\ref{a8}) implies that
\be
A_{0h}=A_{x^h}+\Gamma^i_{lk}y^ky^lA_{ih}.\label{a9}
\ee
Multiplying  (\ref{a9})  with $A^{hj}$, we have
\be
\Gamma^j_{lk}y^ky^l=(A_{0h}-A_{x^h})A^{hj}.\label{a10}
\ee
By  (\ref{YY}) and  (\ref{a10}), one can obtain that
\be
G^j=\frac{1}{2}\Gamma^j_{lk}y^ky^l.\label{a11}
\ee
This means that $F$ is an Antonelli metric.
\end{proof}

\bigskip

\noindent{\bf{Proof of Theorem \ref{main thm 2}}}: By Lemma \ref{Lem1} and Lemma \ref{Lem2}, we get the proof of Theorem \ref{main thm 2}.
\qed

\section{Proof of Theorem \ref{main thm 3}}
\bigskip

\noindent{\bf{Proof of Theorem \ref{main thm 3}}}:
Let $F=\sqrt[m]{A}$ be an $m$-th root Finsler metric and be of isotropic mean Berwald curvature, i.e.,
\be \label{IMB}
{E}_{ij}= \frac{n+1}{2} c F^{-1}h_{ij},
\ee where $c=c(x)$ is a scalar function on $M$. A direct computation implies that the angular metric $h_{ij}=g_{ij}-F^{2}y_iy_j$ are given by the following
\be\label{angular}
h_{ij}= \frac{A^{\frac{2}{m}-2}}{m^2}[mAA_{ij}+(1-m)A_iA_j].
\ee
Plugging (\ref{angular}) into (\ref{IMB}), we get
\be \label{E}
E_{ij}=\frac{(n+1)A^{\frac{1}{m}}}{2m^2 A^2}\ c[mAA_{ij}+(1-m)A_iA_j] .
\ee
By (\ref{YY}), one can see that $E_{ij}$ is rational with respect to $y$. Thus, (\ref{E}) implies that
$c=0$ or
\be\label{E1}
mAA_{ij}+(1-m)A_iA_j=0.
\ee
Suppose that $c\neq0$. Contracting (\ref{E1}) with $A^{jk}$ yields
\[
mA\delta^k_i-A_iy^k=0,
\]
which  implies that $mnA=mA$. This  contradicts our assumption $n\geq 2$. Thus $c=0$ and consequently $E_{ij}=0$.
\qed

\bigskip

\noindent
\textbf{Acknowledgments.} The authors express their sincere thanks to referees and Professor
Z. Shen for their valuable suggestions and comments.

Akbar Tayebi\\
Faculty  of Science, Department of Mathematics\\
Qom University\\
Qom. Iran\\
Email:\ akbar.tayebi@gmail.com
\bigskip

\noindent
Behzad Najafi\\
Faculty  of Science, Department of Mathematics\\
Shahed University\\
Tehran. Iran\\
Email: najafi@shahed.ac.ir

\end{document}